# Is 3D frequency-domain FWI of full-azimuth/long-offset OBN data feasible?

## The Gorgon-data FWI case study

S. Operto[1], P. Amestoy[2], H. Aghamiry[1], S. Beller[1], A. Buttari[3], L. Combe[1], V. Dolean[4], M. Gerest[5], G. Guo[1], P. Jolivet[6], J.-Y. L'Excellent[2], F. Mamfoumbi[1], T. Mary[7], C. Puglisi[2], A. Ribodetti[1], P.-H. Tournier[8]

[1] Université Côte d'Azur, France; email: operto@geoazur.unice.fr,
[2] Mumps Technologies, ENS Lyon, 46 Allée d'Italie, F-69007 Lyon, France; email: patrick.amestoy@mumps-tech.com, jean yves.l.excellent@mumps-tech.com, chiara.puglisi@mumps-tech.com
[3] Université de Toulouse, CNRS, IRIT, Toulouse, France; email: alfredo.buttari@irit.fr
[4] UCA/Univ. of Strathclyde, United Kingdom/UCA; email: Victorita.Dolean@univ-cotedazur.fr
[5] EDF, email: matthieu.gerest@edf.fr
[6] IRIT-CNRS, University of Toulouse, France; email: pierre@joliv.et
[7] Sorbonne Université, CNRS, LIP6, Paris, F-75005, France, email: theo.mary@lip6.fr
[8] Sorbonne University, CNRS, LJLL, France, email: tournier@ljll.math.upmc.fr

Frequency-domain Full Waveform Inversion (FWI) is potentially amenable to efficient processing of full-azimuth long-offset stationary-recording seabed acquisition carried out with a sparse layout of ocean bottom nodes (OBNs) and broadband sources because the inversion can be performed with a few discrete frequencies. However, computing the solution of the forward (boundary-value) problem efficiently in the frequency domain with linear algebra solvers remains a challenge for large computational domains involving tens to hundreds of millions of parameters. We illustrate the feasibility of 3D frequency-domain FWI with a subset of the 2015/16 Gorgon OBN dataset in the NorthWestern shelf, Australia. We solve the forward problem with the massively-parallel multifrontal direct solver MUMPS, which includes four key features to reach high computational efficiency: An efficient parallelism combining message-passing interface and multithreading, block low-rank compression, mixed precision arithmetic and efficient processing of sparse sources. The Gorgon subdataset involves 650 OBNs that are processed as reciprocal sources and 400,000 sources. Mono-parameter FWI for vertical wavespeed is performed in the visco-acoustic VTI approximation with a classical frequency continuation approach proceeding from a starting frequency of 1.7 Hz to a final frequency of 13 Hz. The target covers an area ranging from 260 km$^2$ (frequency ≥ 8.5 Hz) to 705 km$^2$ (frequency ≤ 8.5 Hz) for a maximum depth of 8 km. Compared to the starting model, FWI dramatically improves the reconstruction of the bounding faults of the Gorgon horst at reservoir depths as well as several intra-horst faults and several horizons of the Mungaroo formation down to a depth of 7 km. Seismic modeling reveals a good kinematic agreement between recorded and simulated data but amplitude mismatches between the recorded and simulated reflection from the reservoir suggesting elastic effects. Therefore, future works involve multiparameter reconstruction for density and attenuation before considering elastic FWI from hydrophone and geophone data.

## INTRODUCTION

Full-azimuth long-offset stationary-recording seabed acquisition carried out with a sparse layout of multi-component (4C) ocean bottom nodes (OBNs) and broadband sources are emerging as the leap-forward towards high-resolution imaging in deep water environments by Full Waveform Inversion (FWI). While nowadays FWI is mostly implemented in the time domain, it is well acknowledged that frequency-domain FWI is amenable to efficient processing of long-offset acquisitions because the inversion can be limited to a few discrete frequencies (e.g. Virieux and Operto, 2009). This results because the broad angular illumination of the subsurface provided by such surveys leads to a redundant sampling of the local wavenumber vectors. Here, a local wavenumber vector $\mathbf{k} = (k_x, k_y, k_z)$ denotes the spectral components of a diffractor point that are mapped by one source-receiver pair and one frequency at a given position in the subsurface during FWI (Let's remember that FWI reconstructs the subsurface by summing the images of a fine grid of point diffractors by virtue of the Huygens' principle). This redundant wavenumber sampling provided by different combinations of frequency and source-receiver pair can be reduced by limiting FWI to a few frequencies for the sake of computational efficiency and compact volume of data. The main bottleneck of frequency-domain FWI is generally considered as being related to the forward problem. The frequency-domain forward problem is a boundary-value problem, which requires the solution of a large, sparse, indefinite linear system per frequency, whose solutions are monochromatic wavefields and the right-hand sides (RHSs) are monochromatic sources:

$$\mathbf{AP} = \mathbf{F}, \tag{1}$$

where $\mathbf{A}$ is the so-called impedance matrix resulting from the discretization of the time-harmonic wave-equation operator and the matrices $\mathbf{P}$ and $\mathbf{F}$ gather all the wavefields and sources in their columns. This fundamentally differs from the forward problem of time-domain FWI, which is an initial-value problem mostly tackled with matrix-free explicit time-stepping schemes such as those implemented in reverse time migration (RTM).

The linear system, equation 1, can be solved with two broad categories of linear algebra methods: direct methods (Duff et al., 1986) and iterative methods (Saad, 2003). Direct methods provide solutions in a fixed number of operations. They first perform a RHS-independent pre-processing step (Lower-Upper (LU) decomposition of the matrix, $\mathbf{A} = \mathbf{LU}$) before computing the solutions by forward/backward elimination ($\mathbf{LY} = \mathbf{F}$ followed by $\mathbf{UP} = \mathbf{Y}$). It is often erroneously stated that the memory demand of this pre-processing, related to the storage of the LU factors, prevents realistic 3D wave simulations. This study aims to rectify this shortcut by processinf a real OBN dataset with a leading-edge direct solver. Then, wavefields are computed very efficiently by forward/backward elimination for multiple sources with parallel basic linear algebra subprograms (BLAS), a key feature for FWI. This is a clear advantage over time-domain methods, which require to restart the simulation from scratch each time a new source is processed. The resulting computational burden is classically mitigated by source subsampling or source encoding at the expense of the FWI convergence speed and quality of the reconstructed model (Warner et al., 2013). The second category of linear algebra methods relies on iterative methods, which start from an initial guess of the solution and estimate successive approximations of the solution through a projection process (Saad, 2003). Their main drawback is the convergence speed considering that the matrix **A** is indefinite. Acceptable convergence speed (namely, when the number of iterations scales to frequency) requires efficient preconditioners. We refer the reader to Tournier et al. (2022) where a recent domain-decomposition preconditioner is assessed with realistic 3D benchmarks involving up to a billion unknowns. Compared to direct methods, the advantages of iterative methods are their high parallel efficiency and their moderate memory demand, while their drawback, regardless the convergence issue, is a lower efficiency in processing multiple sources. The pros and cons of direct and iterative methods reviewed above highlight their complementarity to tackle a broad range of 3D FWI applications: Direct solvers may be favored for computational meshes of moderate size and dense acquisitions, while iterative solvers may be more suitable for large targets and sparse acquisitions. Let's conclude this review by recalling other advantages of the frequency domain. The spectral localization provides the most versatile framework to mitigate nonlinearity (cycle skipping) with multiscale strategies based on frequency continuation. Implementation of attenuation in the forward and inverse problems is straightforward and computationally free, while it generates a computational overhead by a factor of two to three in the time domain (Plessix, 2017). Second-order optimization algorithms such as the truncated Newton method are manageable in the frequency domain (Métivier et al., 2013). The source signature can be easily updated at each FWI iteration in alternation with subsurface properties. Finally, subsurface discretization can be matched to the processed frequency according to the resolution power of FWI for sake of computational efficiency and regularization (Operto et al., 2015). However, the frequency domain lacks versatility to implement time-continuation strategies that enable processing early arrivals before reflections through adaptive time windowing. However, time windowing can be replaced by exponential time damping implemented with complex-valued frequencies (e.g. Virieux and Operto, 2009).

In this study, we illustrate the potential of 3D visco-acoustic VTI frequency-domain FWI with a subset of the 2015/2016 Gorgon OBN dataset in the North West shelf, Western Australia. The forward problem is performed with the MUltifrontal Massively Parallel (direct) Solver MUMPS. In the method section, we review the basic principles of frequency-domain FWI and the key functions of the MUMPS solver that provide high computational efficiency. Then, we introduce the Gorgon experiment before proceeding with the FWI results, which reveal the fine structure of the Gorgon horst. Finally, we discuss the computational efficiency of our FWI technology with this dataset.

## METHOD

### Frequency-domain FWI: basic principles

#### *On the computational complexity of time-domain and frequency-domain FWIs*

As above reviewed, 3D frequency-domain FWI should be implemented with a limited number of frequencies to be computationally attractive irrespective of which kinds of solver is used. Theoretical complexity analysis has shown that the number of floating point operations required by either time-domain or frequency-domain FWI applied to dense surface acquisitions ($n^2$ sources) scales to $n^6$ where $n$ denotes one dimension of a $n^3$ subsurface domain. Indeed, the computational time of a wave simulation in the time domain scales to $nt \times n^3$ where $nt$ denotes the number of time steps and hence the time complexity of the forward problem for $n^2$ sources scales to $n^6$ assuming $nt \approx n$. In the frequency domain, the computational time of the LU factorization of the matrix **A** scales to $n^6$, its memory storage scales to $n^4$ and hence the time complexity of the solution step for $n^2$ sources scales to $n^4 \times n^2 = n^6$ (the factorization and the solution steps have the same time complexity for $n^2$ RHSs). For iterative methods, the computational time of a sparse matrix-vector product scales to $n^3$, the number of iterations generally scales to $n$ leading to a $n^6$ time complexity for $n^2$ sources. This prompts us to conclude that the relative efficiency of time-domain versus frequency-domain FWI cannot be drawn from this basic complexity analysis. Only numerical experiments carried out with massively-parallel solvers and realistic benchmarks are useful to choose the most suitable FWI implementation according to the experimental setting as defined by the size of the computational domain, the number of RHSs, the acquisition geometry, the frequency bandwidth involved in FWI and the physics of wave propagation.

*On the frequency management in frequency-domain FWI*

In frequency-domain FWI, the discrete frequencies are typically processed hierarchically from the low frequencies to the higher ones to design multiscale imaging and mitigate nonlinearity resulting from cycle skipping. Ideally, one can process a narrow batch of frequencies rather than a single frequency during a multiscale step to preserve wavenumber redundancy during multi-parameter FWI and mitigate parameter cross-talk accordingly, or improve signal-to-noise ratio by increasing fold. In this setting, several frequencies can be processed efficiently by distributing them over groups of processors. An analog parallelism is classically used in time-domain FWI by distributing a batch of sources over group of processors. Typically, the group of processors manages the domain decomposition of the subsurface in the time domain, while it manages the partitioning of the graph of the matrix in the frequency domain. Also, several cycles of frequency-domain FWI can be performed to reprocess the low frequencies at cycle *n +1* starting from the final subsurface model of cycle *n*. This outer loop is a last resort to balance the fact that the full frequency band cannot be processed in one go during the last multiscale step of frequency-domain FWI unlike in time-domain FWI.

## Forward problem: On the discretization of the wave equation

When the linear system, equation 1, is solved with a direct solver, the LU decomposition of the sparse impedance matrix **A** generates some significant fill-in (namely, additional non-zero coefficients are generated), which is minimized by fill-reducing matrix ordering techniques such as nested dissection methods. This re-ordering should be however combined with a discretization method that maximizes the sparsity of the matrix and minimize its numerical bandwidth while providing a sufficient accuracy for a discretization rule of four grid points per wavelength (the coarsest parametrization sampling half the wavelength, i.e., the size of the smallest structure that can be captured by FWI). We satisfy these specifications with the finite-difference method of Operto et al. (2014) recently improved by Aghamiry et al. (2022). Operto et al. (2014) decompose the VTI acoustic wave equation as the sum of an elliptic anisotropic wave equation and an anelliptic term for horizontal pressure while the vertical pressure has a closed-form expression as a function of the horizontal pressure. This allows us to perform seismic modeling in VTI media without computational overhead compared to the isotropic case. Moreover, we design a compact stencil with a consistent mass matrix and a weighted combination of stiffness matrices discretized with several second-order accurate stencils. Finally, a high uniform accuracy of the wavefields is provided by adaptive weights minimizing numerical dispersion at each point of the subsurface according to the local wavelength (Aghamiry et al., 2022).

## Forward problem: Pushing the limits with the MUMPS solver

Once a suitable discretization scheme has been designed, the computational cost of frequency-domain FWI is, by and large, controlled by the computational efficiency of the solver that is used to solve the multi-RHS system, equation 1. The MUMPS solver contains four key ingredients that allow us to tackle large-scale FWI case studies: multi-level parallelism combining Message Passing Interface (MPI) and multithreading, Block Low-Rank (BLR) approximation, mixed precision (MP) arithmetic and efficient processing of multiple sparse RHSs.

### *The massively-parallel multifrontal method*

The multifrontal method (Duff et al., 2017) recasts the LU factorization of the sparse matrix **A** as a sequence of dense factorizations of the so-called frontal matrices. This sequence is defined by a tree-shaped dependency graph called the elimination tree, which has a front associated with each of its nodes. In a parallel environment, the workload is distributed onto the processors according to a mapping of the processors on the nodes of the elimination tree (Amestoy et al., 2016, their figure 1). Two sources of parallelism, referred to as tree parallelism and node parallelism, are then exploited by the MPI processes and at a finer level with multithreading.

### *Block Low-Rank (BLR) multifrontal solver*

Bebendorf (2004) has shown that matrices resulting from the discretization of elliptic partial differential equations, such as the Helmholtz equation, have low-rank properties. In the context of the multifrontal method, frontal matrices are not low-rank themselves but exhibit many low-rank sub-blocks. This low-rank property can be used to reduce the storage and the number of floating-point operations both during the LU and solution steps. In order to achieve a satisfactory reduction in both the computational complexity and the memory footprint, sub-blocks have to be chosen to be as low-rank as possible. This can be achieved by clustering the unknowns in such a way that an admissibility condition is satisfied (Amestoy et al., 2016, their figure 4). This condition states that a sub-block interconnecting different variables will have a low rank if these associated variables are far away in the domain, because they are likely to have a weak interaction. These blocks are compressed with a truncated QR factorization with column pivoting using the user-defined threshold parameter $\varepsilon$. Typically, $\varepsilon$ ranging between $10^{-4}$ and $10^{-5}$ can be used for FWI. Thanks to the low-rank compression, the theoretical complexity of the factorization is reduced from $O(n^6)$ to at best $O(n^5)$ and the size of the LU factors is reduced from $O(n^4)$ to $O(n^{3.5})$ (Amestoy et al., 2018b).

*Mixed precision arithmetic Block Low-Rank (BLR) method*

The recent and growing support of low precision arithmetic on modern hardware, such as half precision floating-point arithmetic, provides new opportunities for accelerating direct methods and reducing their memory consumption. The Block Low-Rank (BLR) approach described above is particularly amenable to mixed precision arithmetic, because it allows for easily splitting each block into multiple precisions based on its singular values (Amestoy et al., 2022). In the current version of MUMPS, this mixed precision representation exploits 16-bit and 24-bit formats to reduce memory consumption only; computations are still performed entirely in single (32-bit) precision arithmetic. However, the volume of communications is also reduced, which can positively affect the time performance. Moreover, with mixed precision BLR activated, the low-rank admissibility condition is modified so that more blocks can be represented under low-rank form; this decreases the number of floating point of operations and can also accelerate the computation.

*Exploiting the sparsity of the seismic sources*

Once the LU factorization has been performed, direct solvers efficiently process blocks of RHSs with parallel BLAS. Moreover, the source matrix **F**, equation 1, is extremely sparse for seismic sources. Amestoy et al. (2018a) exploited this sparsity to greatly reduce the number of operations performed and the volume of data exchanged during the solution phase. During the forward phase, a property established by Gilbert and Liu (1993) can be used to limit the number of operations for computing the solution **Y** of **LY** = **F**: for each column $f_j$ of **F**, one has only to follow a union of paths in the elimination tree, each path being defined by a non-zero entry in $\mathbf{f}_j$. The flop reduction results from this property and from the fact that the non-zero entries in each source $\mathbf{f}_j$ have a good locality in the elimination tree so that only a few branches of the tree per source need to be traversed. This may be not enough to provide large gains in a parallel context since processing independent branches of the elimination tree is necessary to provide significant parallelism. This is achieved by permuting the columns of the matrix **F** to reduce the number of operations while providing parallelism. In the framework of FWI, we exploit the fact that the non-zero entries of one block of columns of **F** are well clustered in the elimination tree when the corresponding sources are contiguous in the computational mesh (i.e., follow the natural ordering of the acquisition) and when nested dissection is used for reordering the unknowns of the problem. During the application hereby presented, we exploit this geometrical property of the source distribution to choose those relevant subsets of sparse *RHSs* that are simultaneously processed in parallel during the forward elimination phase.

# Gorgon FWI case study

## Geological context, acquisition and data anatomy

### *Targeted area*

The Gorgon gas field is a giant field located at the southern end of Rankin platform in the Northern Carnarvon basin, North-West Australia, approximately 130 km offshore Barrow Island (Figure **1**a). The geological structure of the Gorgon gas field is an elongated horst block trending north-north east at the south western end of the Rankin trend. The target depth is approximately 3600-4300 m subsea and 2.75-3.3 seconds two-way traveltimes. The horst is bounded by steep faults on the eastern flank with throws of up to 1.5 km, and similar faults, with slightly smaller throws on the western flank. Several intra-horst en-echelon faults, with throws of tens to hundreds of meters, separate the field into fault blocks (Figure **1**a, top-left inset). The core area is more elevated and is referred to as the main horst (fault blocks 1 to 5 in Figure **1**a).

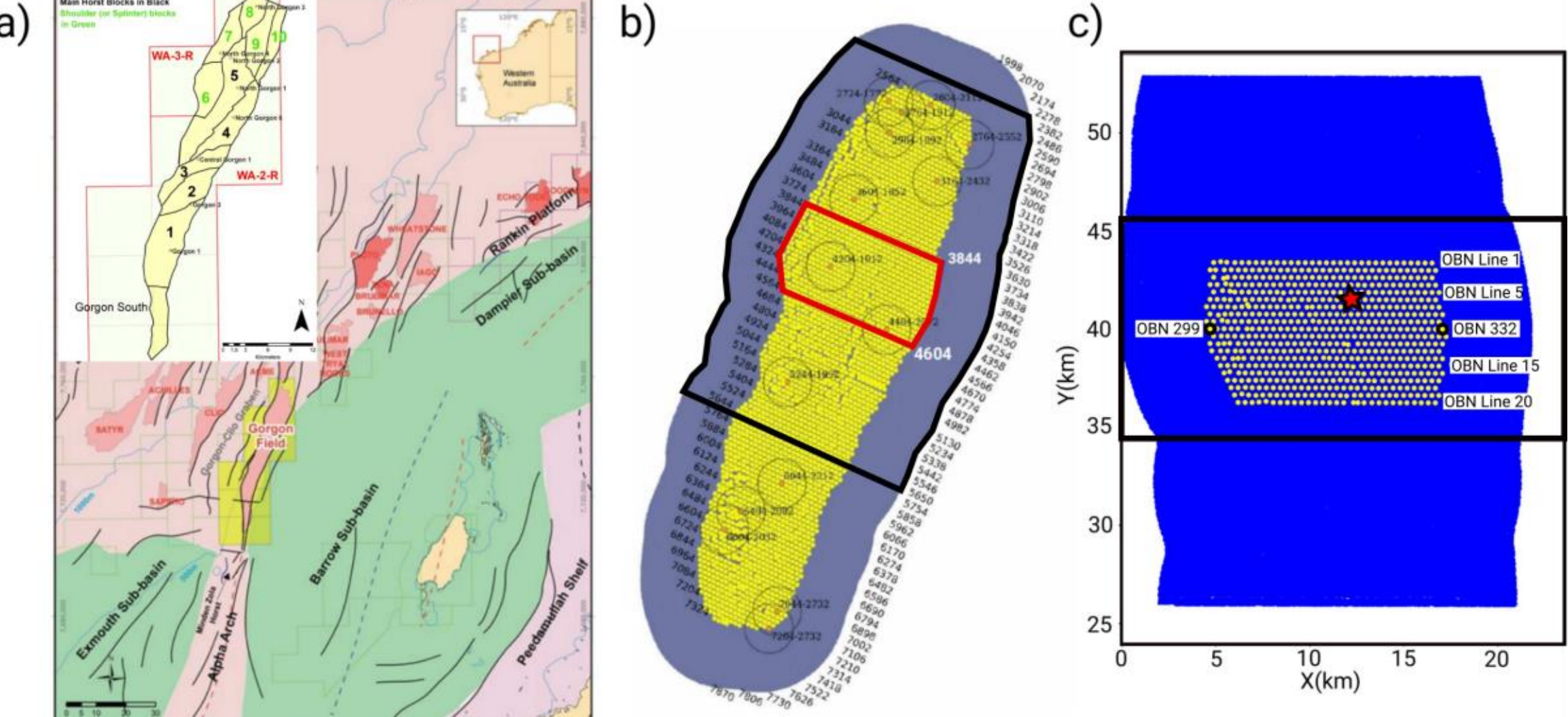

Figure 1: a) Structural diagram of targeted area (Gorgon field). Top-left inset shows fault blocks partitioning Gorgon horst (courtesy M Van Ross, Chevron). b) Shot carpet and OBN layout of 2015/2016 Gorgon survey. Thick red and black polygons delineate the subset of OBNs and shot used in this study. (c) FWI target after rotation/translation of coordinate system for frequencies ≤ 8.5 Hz (target 1). OBN are plotted with yellow circles. Shot carpet is plotted in blue. OBNs 299 and 332 analyzed in this study are highlighted. Star shows well position (Figure 8). Thick black rectangle delineates the FWI target for frequencies ≥ 8.5 Hz (target 2). All labels in next figures are provided in this coordinate system.

*The 3D Gorgon OBN survey*

The 3D Gorgon OBN (Gorgon Ocean Bottom Node Seismic Survey) data was acquired in 2015/16 to overcome fault shadows observed from the processing of narrow-azimuth towed-streamer data and to serve as a 4D seismic baseline (Figure.**1**b). The water depth ranges from 100 m to 900 m. The shot carpet, which involves 735 source lines and 697,345 shot points, is approximately 980 $km^2$ with a full imaging area of approximately 240 $km^2$. A total of 3100 4-C OBNs distributed along 120 lines covering an area of 436 $km^2$ were deployed on a staggered 375 m grid. The survey was acquired with 18.75 m × 37.5 m shot interval/source line interval on a dual source flip-flop setting. Twenty OBN lines involving 650 nodes and a source carpet of 400,258 shots were provided to us for FWI application (Figure **1**b). The acquisition map associated with this sub dataset after rotation and translation of the coordinate system is shown in Figure **1**c.

*Initial subsurface model and data anatomy*

A common-OBN gather is shown in Figure 2. Four diving-wave arrivals (d1, d2, d3, d4) and three reflections at critical distances (r2, r3, r4) can be easily identified. A shadow zone in the data beyond 8 km offset highlights a low-velocity zone above the reservoir. The arrivals r4 and d4 denote the reflection from the top of the reservoir and the diving wave from beneath. The critical distance where r4 and d4 join is around 9 km. A VTI subsurface model (vertical wavespeed $V_0$, Thomsen's parameters $\delta$ and ε), discretized on a 75 $m$ × 75 $m$ × 3 $m$ grid was provided to us for FWI (see Figure 3a for the vertical section of $V_0$ across the OBN position shown in Figure 2). Ray tracing from the OBN position shows that arrivals d1, d2, d3 propagate above the low velocity zone with a maximum penetration depth of 1.5km and a maximum offset of around 7 km while the diving wave from the reservoir reach the surface at a minimum offset of 9 km (Figure 3b). These observations are consistent with the anatomy of the data (Figure 2). Arrival PLM in Figure 2 undergoes a reflection from free surface before being reflected by a discontinuity in the overburden (Figure 3b).

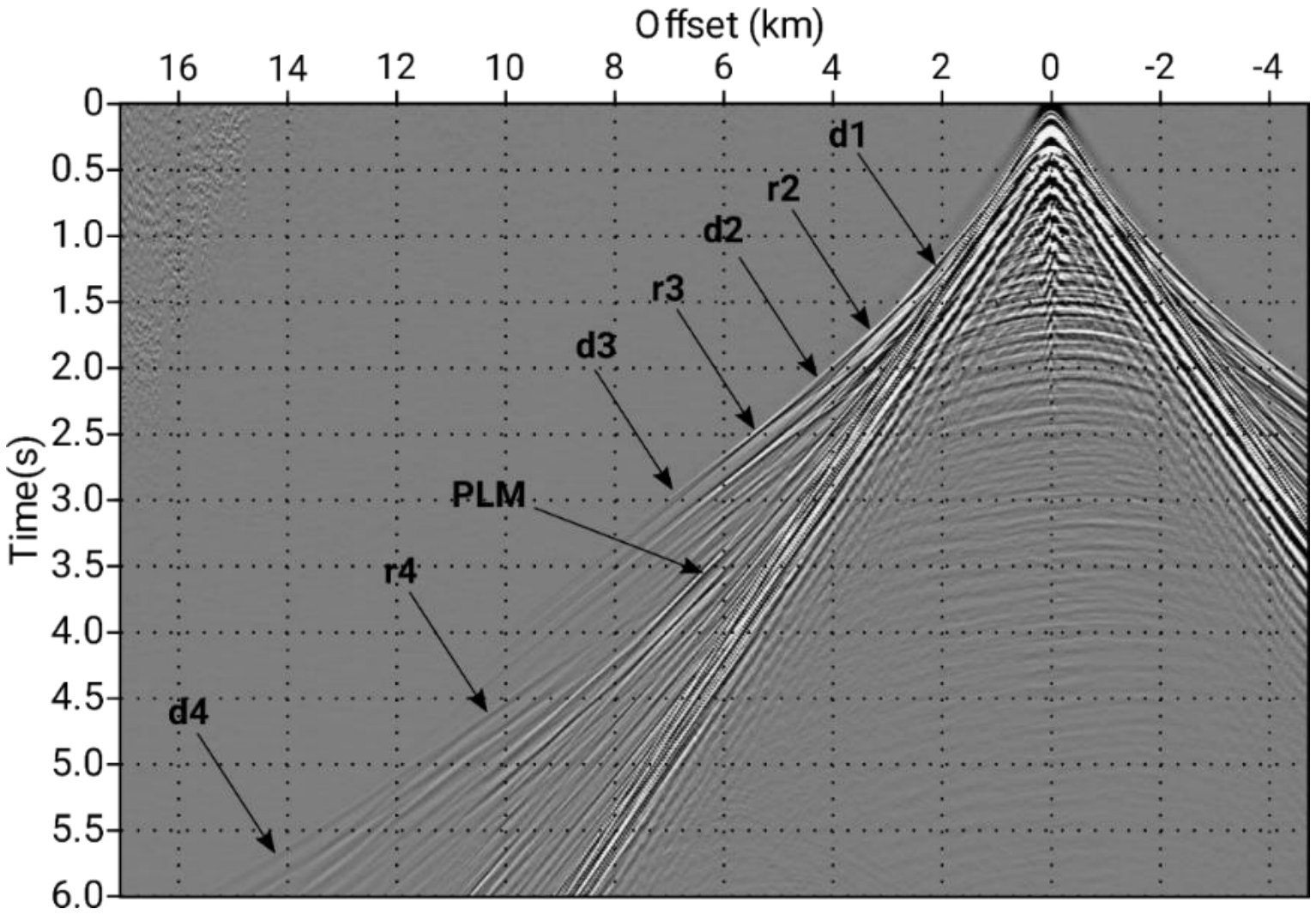

Figure 2: Common-OBN gather for source line intersecting OBN position (OBN 332 in Figure1c). Main diving-wave arrivals (d1, d2, d3, d4) and reflections at critical distances (r2, r3, r4) are pointed by arrows. Arrival d4 is reflection from top of reservoir. Arrival PLM is illustrated in Figure 3b.

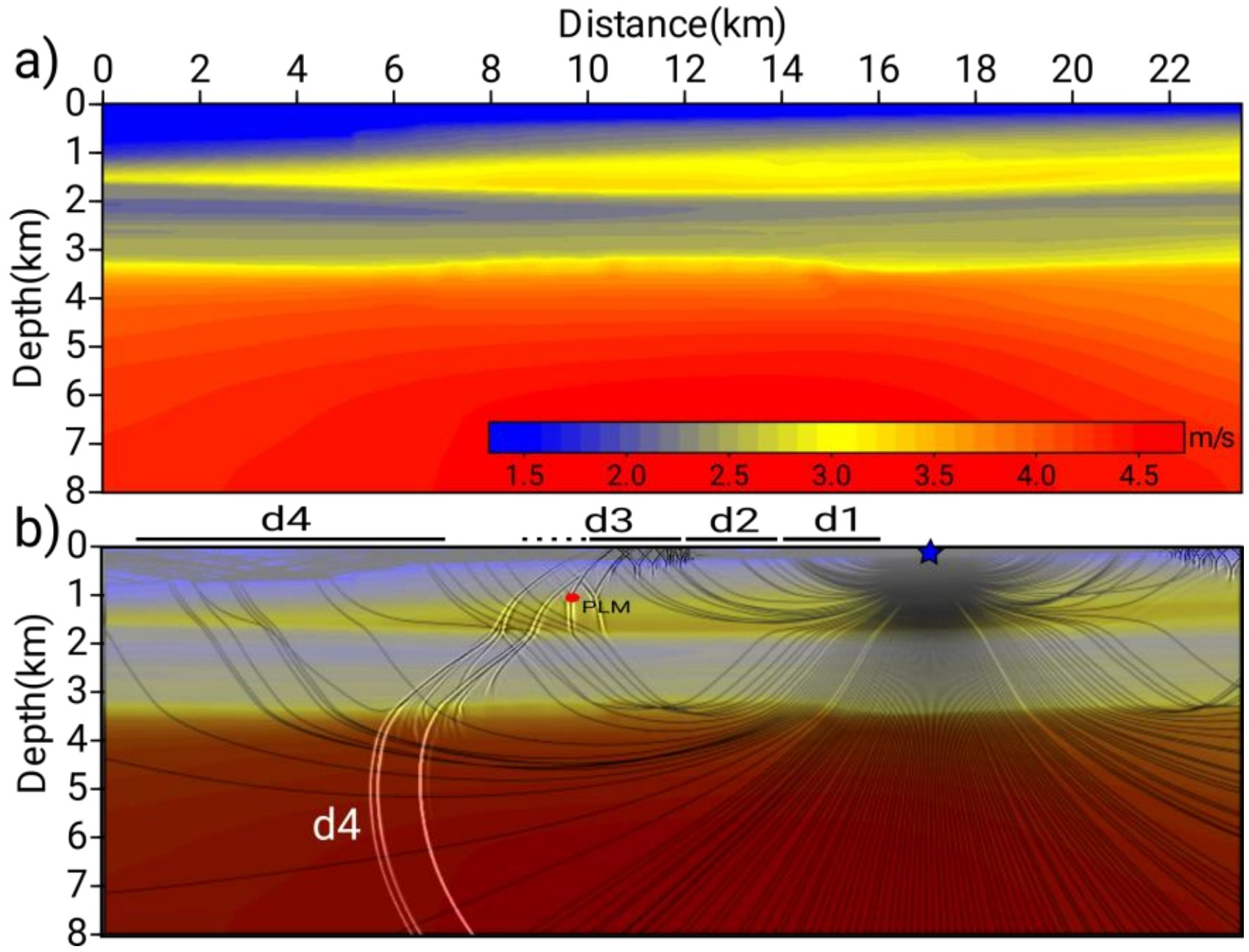

Figure 3: a) Dip section of starting model ($V_0$) along OBN line 10 (Figure 1c). b) Isotropic ray tracing in $V_0$ model with superimposed wavefield snapshot (T=3.6s). The red circle points the reflection point of the PLM arrival. This arrival undergoes a first reflection from the sea surface before a second reflection from a velocity discontinuity at around 1 km depth. The first arrival is an evanescent wave grazing on the top of the low velocity zone and is not visible in the figure. The “apparent” first arrival is the diving wave from the reservoir (d4) at the critical distance (which is underestimated (∼7 km) compared to the observed one (∼10 km) in Figure 2, black arrow r4). The blue star points the source position. The segments labeled by d1, d2, d3, d4 indicate the offset range over which the corresponding diving waves are recorded as a first arrival in Figure 2.

## Frequency-domain FWI setup

We perform successive monoparameter monochromatic inversions from 1.7 Hz to 13 Hz. We update $V_0$ from the pressure data while we process density, attenuation, $\delta$ and $\epsilon$ as passive parameters. We use a constant quality factor of 200 below the sea bottom to build the background attenuation model with the Kolsky-Futterman relationship. The density model is inferred from $V_0$ using the polynomial relation of Brocher (2005). Between 1.7 Hz and 8.55 Hz (13 processed frequencies), we involve all the available short carpet in FWI although a significant part of the shot carpet is not overlaid by the available OBN layout (Figure 1c). First, this allows us to illustrate the size of the targets that can be routinely processed at these low frequencies by frequency-domain FWI and direct solvers. Second, considering the good accuracy of the starting model, the long-offset data improve the illumination of the central area overlaid by the OBNs through undershooting. Between 8.55 Hz and 13 Hz (five processed frequencies), we limit the size of the target to the area delineated by the OBN layout due to restricted access to computational resources (Figure **1**c). The dimensions of the first target are 23.5 $km(dip) \times 30\ km(cross) \times 8\ km(depth)$ while the cross direction is limited to 11 km in the second target. The 650 OBNs, processed as reciprocal sources, are involved at each FWI iteration, taking advantage of the computational efficiency of the solution step of direct solvers. During each mono-frequency inversion, the grid interval is matched to the frequency to satisfy a discretization rule of around four grid point per minimum wavelength. Accordingly, we down-sample the available $V_0$, $\delta$ and $\epsilon$ models on a 150 m grid to start FWI at the 1.7-Hz frequency. To build baseline results and mitigate acquisition footprint more efficiently, we discretize frequencies with small steps (0.5-1 Hz) although the frequency interval could be increased with frequency for sake of computational efficiency according to the wavenumber-coverage rule of Sirgue and Pratt (2004). Free-surface boundary condition is applied on top of the grid, that is surface multiples are involved in FWI, and sources and receivers are positioned at arbitrary positions in the finite-difference grids with Kaiser-windowed sinc functions (Hicks, 2002). We perform the inversion with the quasi-Newton l-BFGS algorithm and the line search of the SEISCOPE optimization toolbox (Métivier and Brossier, 2016). Monochromatic source signatures are estimated in alternating mode with the model perturbations at each FWI iteration. We didn't apply regularization although it may be beneficial to mitigate acquisition footprint. We just applied total-variation (TV) denoising on the final FWI model. The reader is referred to Operto et al. (2015) and Amestoy et al. (2016) for an application of this basic workflow on OBC data from the North Sea while a more elaborated multiparameter FWI for $V_0$, density and attenuation is applied in Operto and Miniussi (2018).

## FWI results and quality control

### *FWI results*

We present our imaging results by showing either the raw FWI model, the magnitude of the velocity gradient, namely $[(dV_0(x,z)/dx)^2 + (dV_0(x,z)/dz)]^{1/2}$, as a pseudo-reflectivity of the subsurface or both by overlaying the pseudo-reflectivity on the FWI model. This later representation is helpful to check the geometrical consistency between the pseudo-reflectivity (faults, stratigraphic horizons) highlighted by the differential information embedded in the FWI model and the main structural units associated with large-scale velocity variations. We first show depth, dip and cross sections of the 23.5 $km$×30 $km$×8 $km$ FWI model (target 1) inferred from the 8.55 Hz inversion (Figure 4). The depth slice clearly shows the steep faults bounding the horst as well as several short-scale structures in the horst. At this frequency, acquisition footprint is quite strong in the depth slice. The dip section clearly shows the overall geometry of the Gorgon horst, its bounding faults and different from the Mungaroo formation. The western flank of the horst is also highlighted in the cross section at the intersection with the dip section.

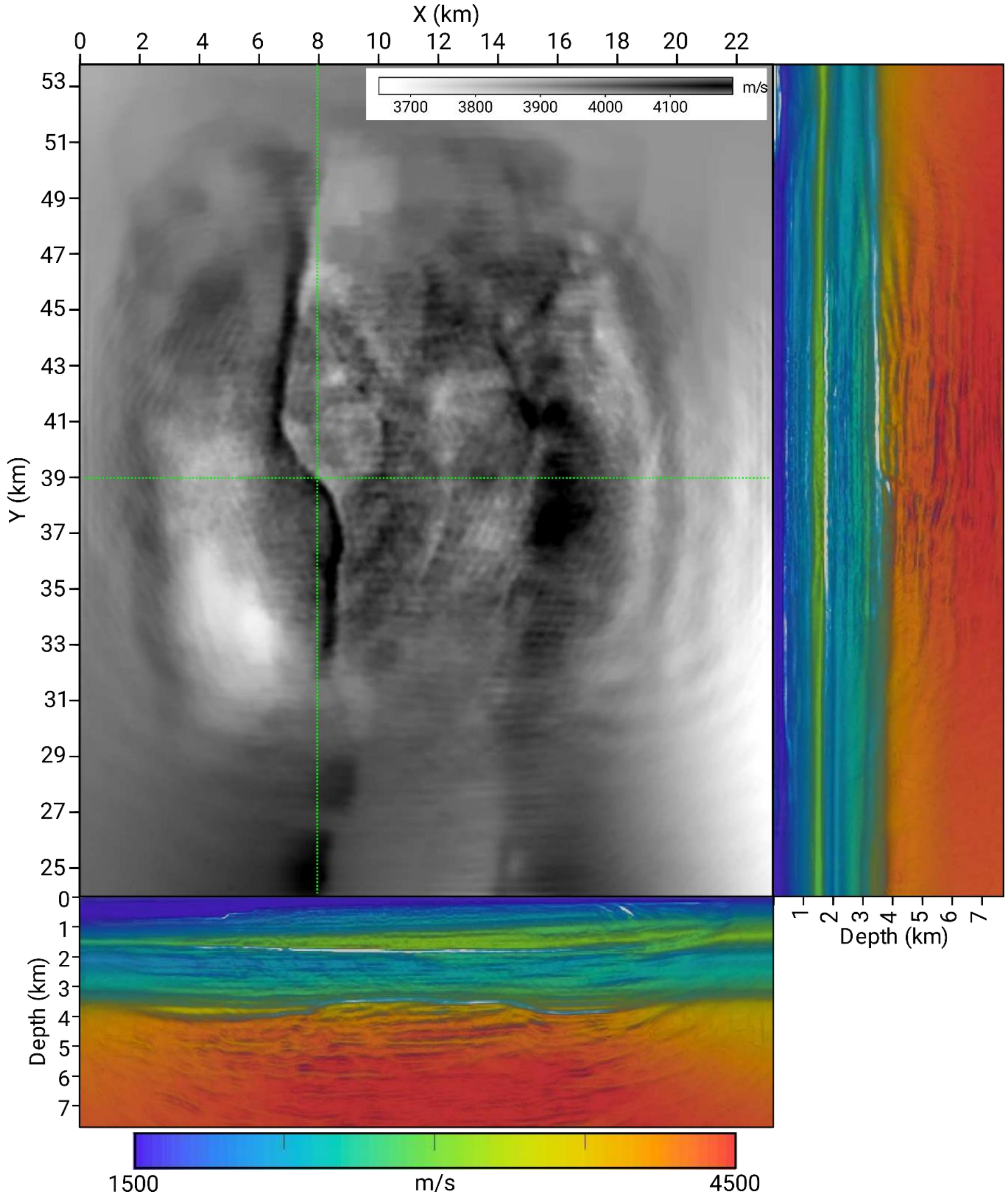

Figure 4: Sections across 8.55-Hz FWI model. Depth slice at 4 km depth is shown with a gray scale. A dip and cross vertical sections are plot with a color scale with magnitude of velocity gradient overlaid. Positions of these sections (Y=39 km, X=8 km) are shown by the dot lines overlaid on depth slice. They cross on the western flank of the horst.

We continue FWI up to 13 Hz to boost spatial resolution and mitigate acquisition footprint in the target 2 (Figure 1c). Indeed, adding the contribution of a wider range of frequencies contributes to better focus not only the subsurface model but also the spurious image of the OBN grid on the seabed, hence mitigating its footprint in the FWI model. Comparing depth slices across the reservoir of the starting and final FWI models highlights the dramatic resolution improvement provided by FWI (Figure **5**).

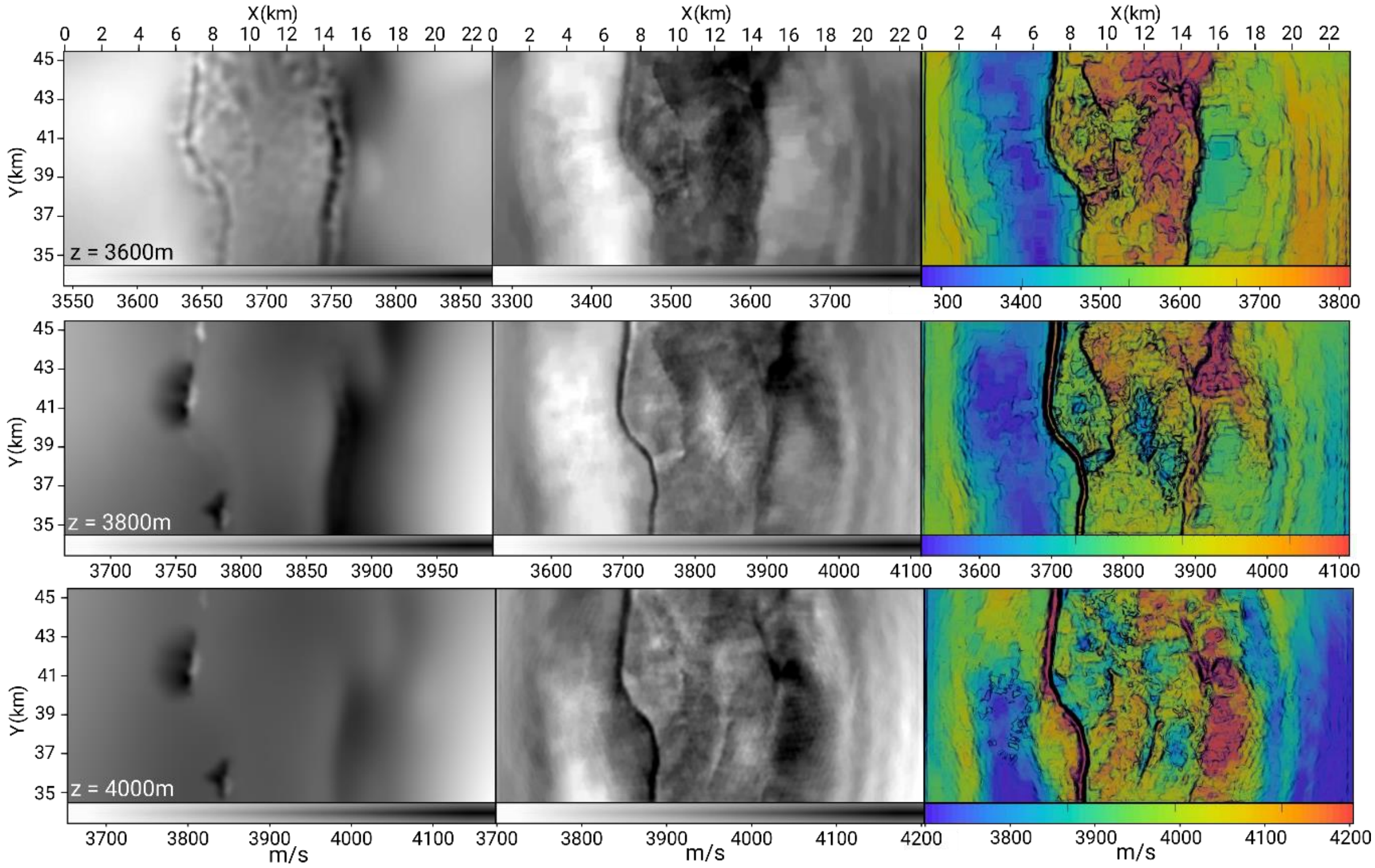

Figure 5: Depth slices at reservoir level. From top to bottom, depths are 3600 m, 3800 m and 4000 m. Left and middle columns show slices of starting and final 13-Hz FWI models ($V_0$) with a gray scale. Right column shows slices of final FWI model with color scale. Magnitude of velocity gradient is overlaid to highlight contrasts.

Compared to Figure **4**, the bounding faults of the horst have been sharpened and the acquisition footprint has been mitigated. Moreover, intra-horst en-echelon faults can now be delineated quite clearly. In particular, the geometry of the intra-horst fault separating the blocks 5 and 6 (Figure **1**a, inset) is well delineated. A dip section of the starting and final FWI models can be compared in Figure **6**(a-b). We show separately the pseudo-reflectivity built by the magnitude of the velocity gradient and the sum of the vertical and horizontal velocity derivatives, namely $dV_0(x,z)/dx + dV_0(x,z)/dz$, in Figure **6**(c-d) to better see the geometry of the intra-horst faults and the stratigraphic horizons between 3600 m and 7000 m depth. The two intra-horst faults separating the main block from the blocks 4 and 6 are well identified (Figure **1**a, inset). Two other dip sections of the FWI model are shown in Figures **7** to highlight structural variations in the cross direction. Comparison of the FWI model with the North-Gorgon 1 well log (Figure **1**c) confirms that FWI manages to capture several horizons in the reservoir (Figure 8).

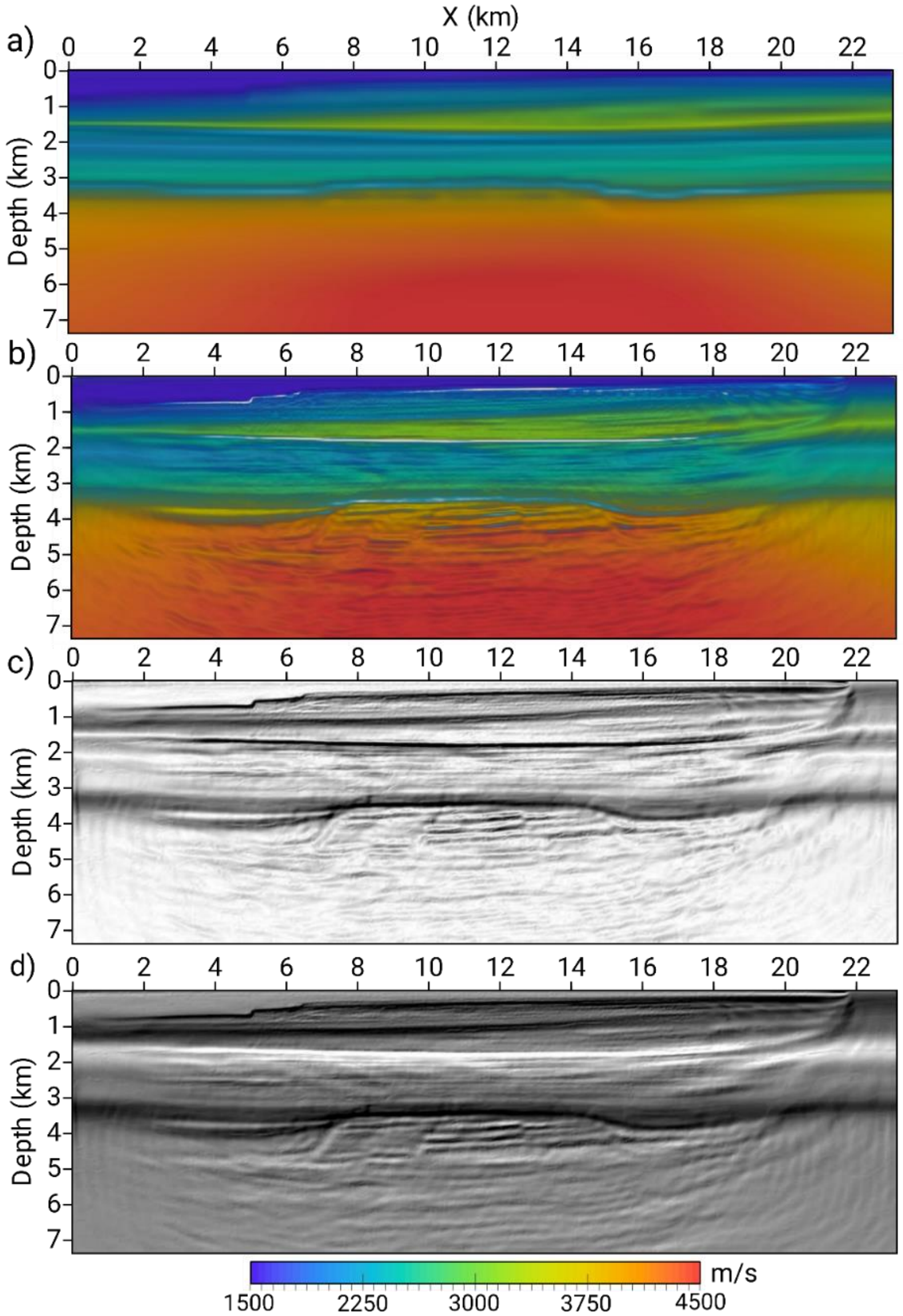

Figure 6: (a) Dip section of starting $V_0$ model (Y = 40 km)). (b) Section of final FWI model with magnitude of velocity gradient superimposed. (c) Magnitude of velocity gradient of final FWI model. (d) Sum of vertical and horizontal derivatives of final FWI model.

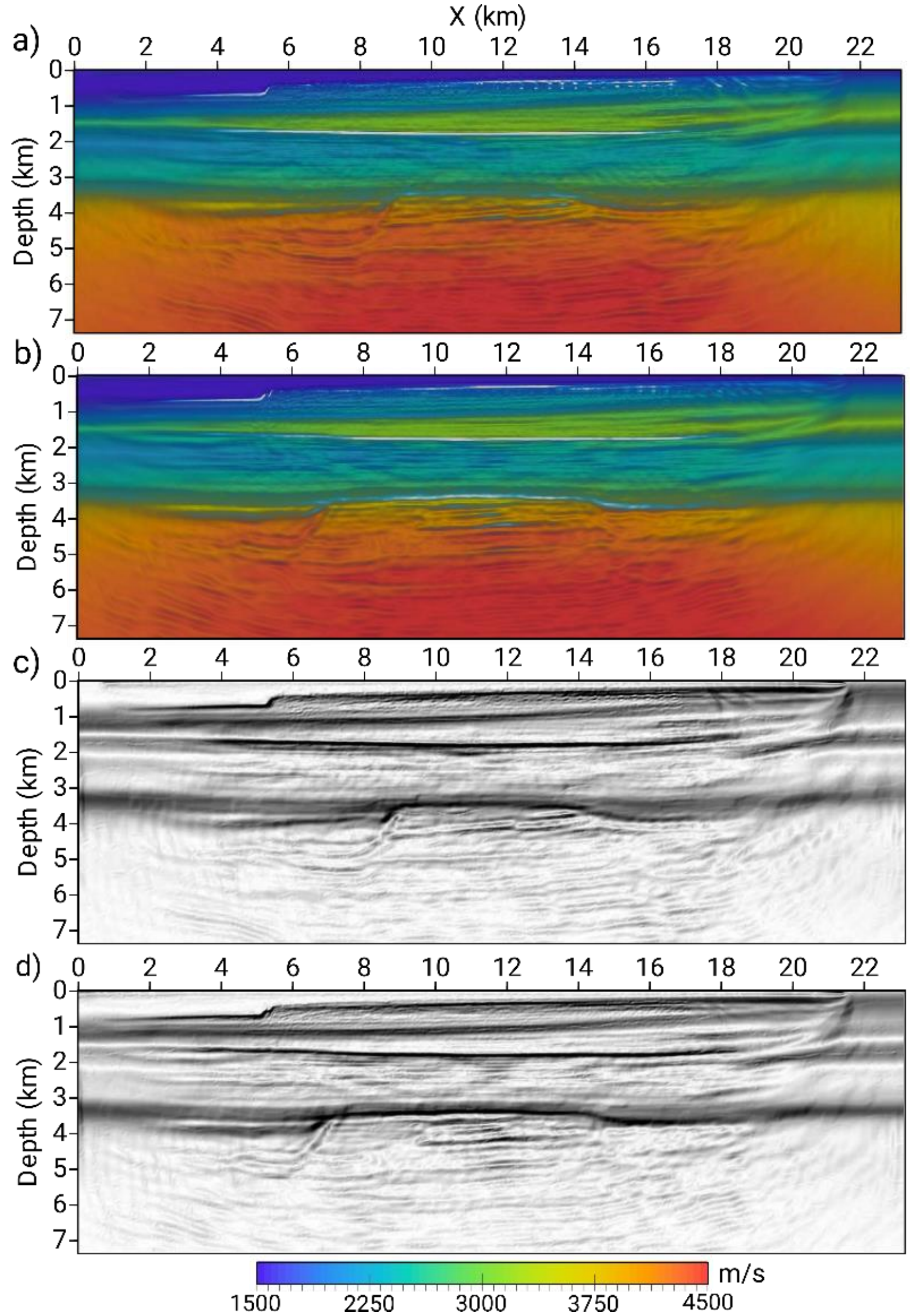

Figure 7: (a-b) Dip sections of final FWI model with magnitude of velocity gradient overlaid. (a) Y=37 km. (b) Y=42 km. (c-d) Magnitude of gradient of final FWI model for sections shown in (a-b).

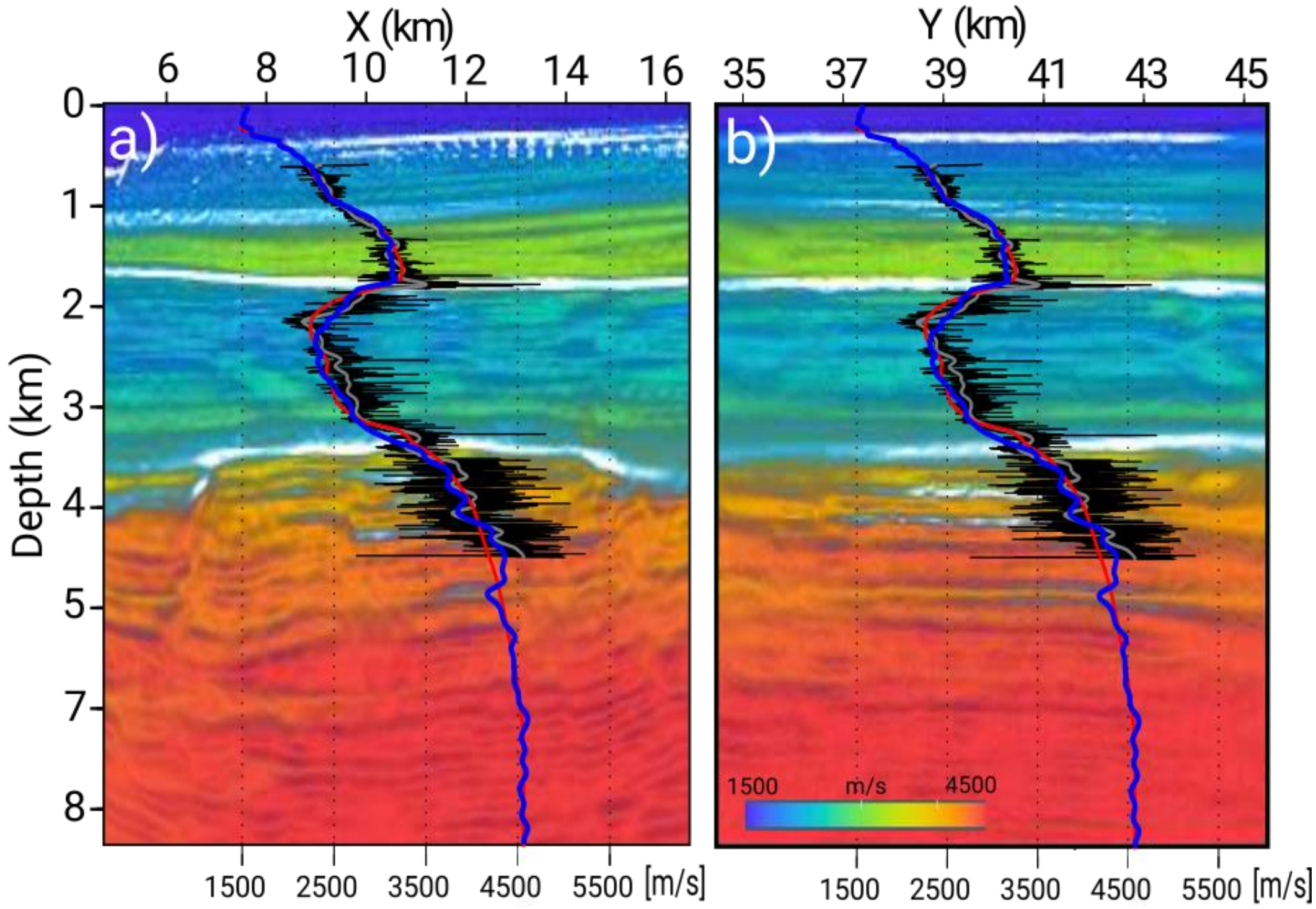

Figure 8: (a) Part of the dip section of the final FWI model with superimposed well log (dot lines). (b) Cross section with superimposed well log. Well log, its smooth version, and profiles of starting and FWI model are plot in black, green, red and blue, respectively. Note reconstructed perturbations at reservoir level (beneath 3.6 km depth).

*Quality control*

We first assess data fit in the frequency domain. Comparison between recorded and simulated data at the 1.7 Hz starting frequency is shown in Figure **9** for OBN 299 (Figure 1c). Although the data are noisy at this frequency, we were able to perform few iterations to update the shallow part of the structure. Another example is shown in Figure **10** for the 5-Hz frequency. In the dip (X) direction, the match of the amplitudes and phases is acceptable up to 10 km offset. Beyond this offset, the simulated amplitudes are systematically overestimated. This offset corresponds to the critical distance of the reflection from the top of the reservoir (Figure **2**). This amplitude mismatch might result from elastic effects.

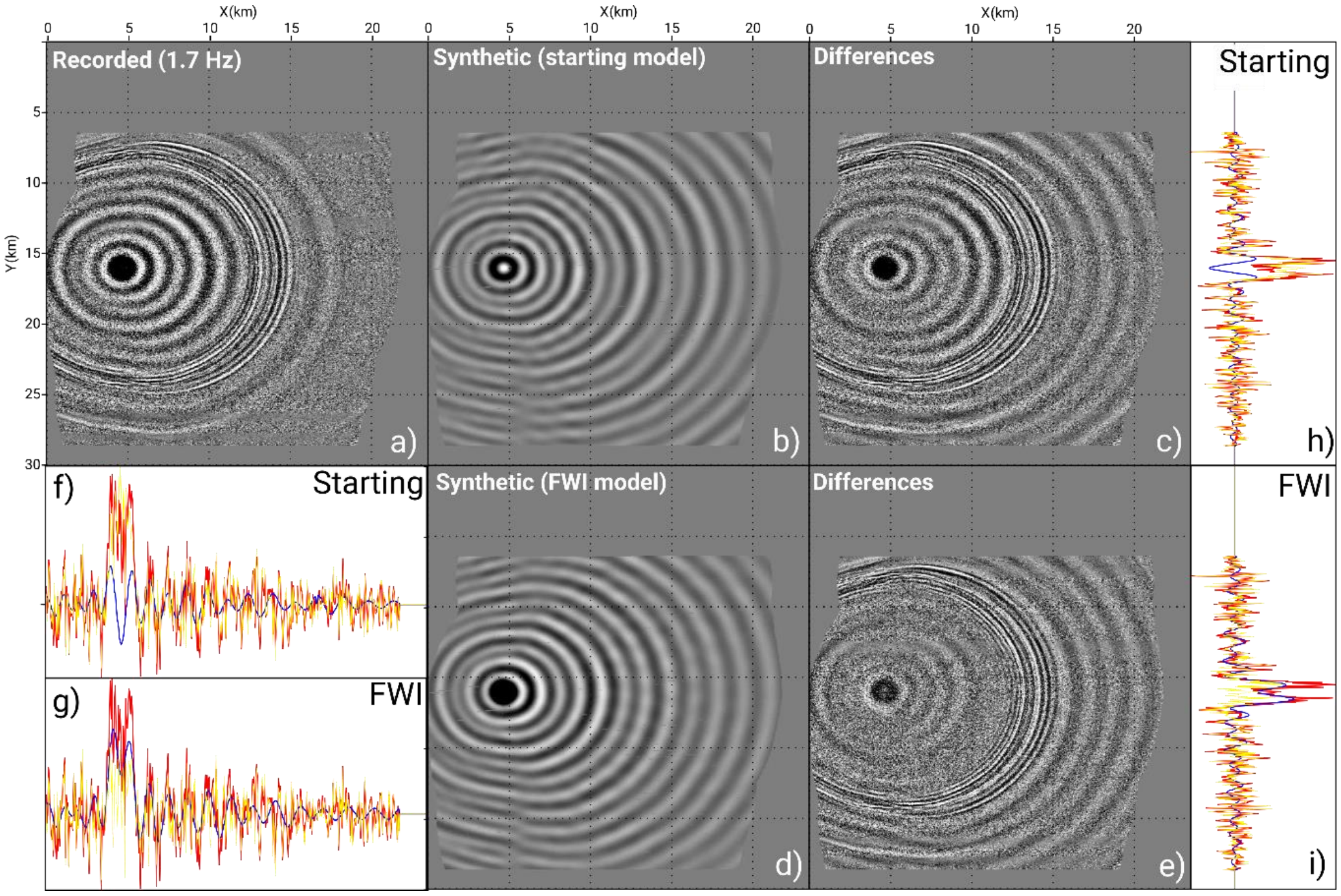

Figure 9: Frequency-domain common-OBN gather 299 (real part). Frequency is 1.7 Hz. a) Recorded data. b) Synthetics computed in starting model. c) Differences. (d-e) Same as (b-c) but synthetics are computed in 1.7-Hz FWI model. (f-i) Direct comparison between recorded data (red) and computed data (blue), and their differences (yellow). (f-g) X profile. (h-i) Y profile. Synthetics are computed in starting (f,h) and final FWI (g,i) model.

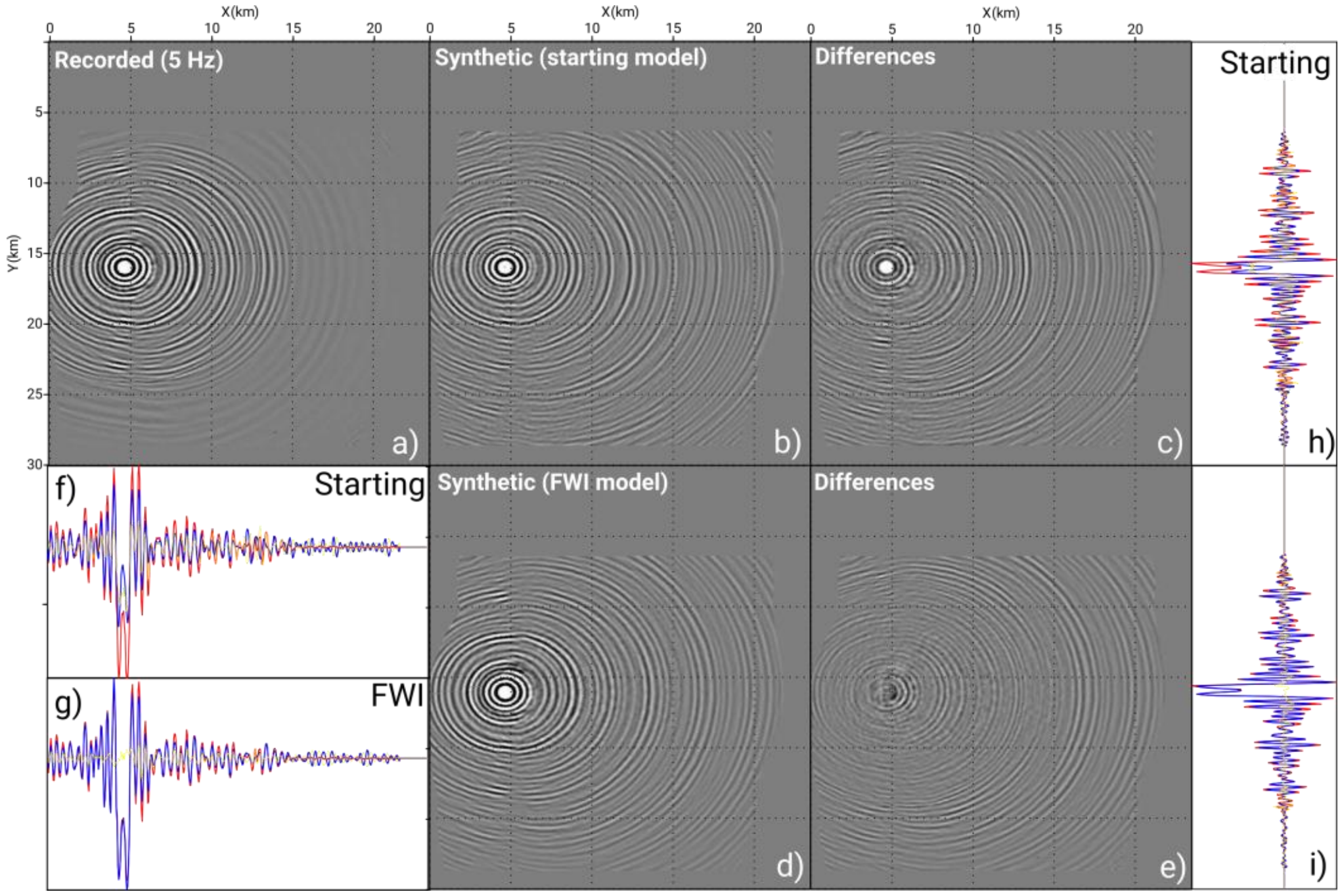

Figure 10: Same as Fig 9 for 5-Hz frequency.

Recorded and simulated time-domain common OBN gathers are shown in Figures **11** for OBN 299 (Figure **1**c). The recorded and synthetic data are band-pass filtered in the 2-14 Hz frequency band. The simulated data are shown after convolution with the estimated source signature by matching the recorded data. Synthetics are computed in the starting and final FWI model to assess the data-fit improvement achieved by FWI. In the bottom row of Figures **11**, we interleave 20 consecutive recorded traces with the simulated counterparts to assess the kinematic agreement between the two sets of seismograms and identify which arrivals have been reconstructed by FWI. We plot these seismograms with a linear move out to compress the time scale and highlight traveltime mismatches. First-order conclusions are that traveltimes of early arrivals and reflection from the reservoir are well matched after FWI. Amplitudes of the simulated post-critical reflection from the top of the reservoir are overestimated suggesting elastic effects as already suggested by frequency-domain modeling. However, we note that the critical distance and the AVO trend of this reflection are better matched after FWI (Figure **11**, arrows). "Conversely, amplitudes of the simulated short-spread reflections from the overburden are underestimated. The match of the well logs in Figure 8 shows that the contrast at 1.7-km depth has been softened by FWI resulting in the above-mentioned underestimated reflection amplitudes. This softening might result in turn from the inaccurate background density and/or attenuation models or from slight mispositioning in depth of the contrast already present in the initial velocity model at 1.7 km depth.”

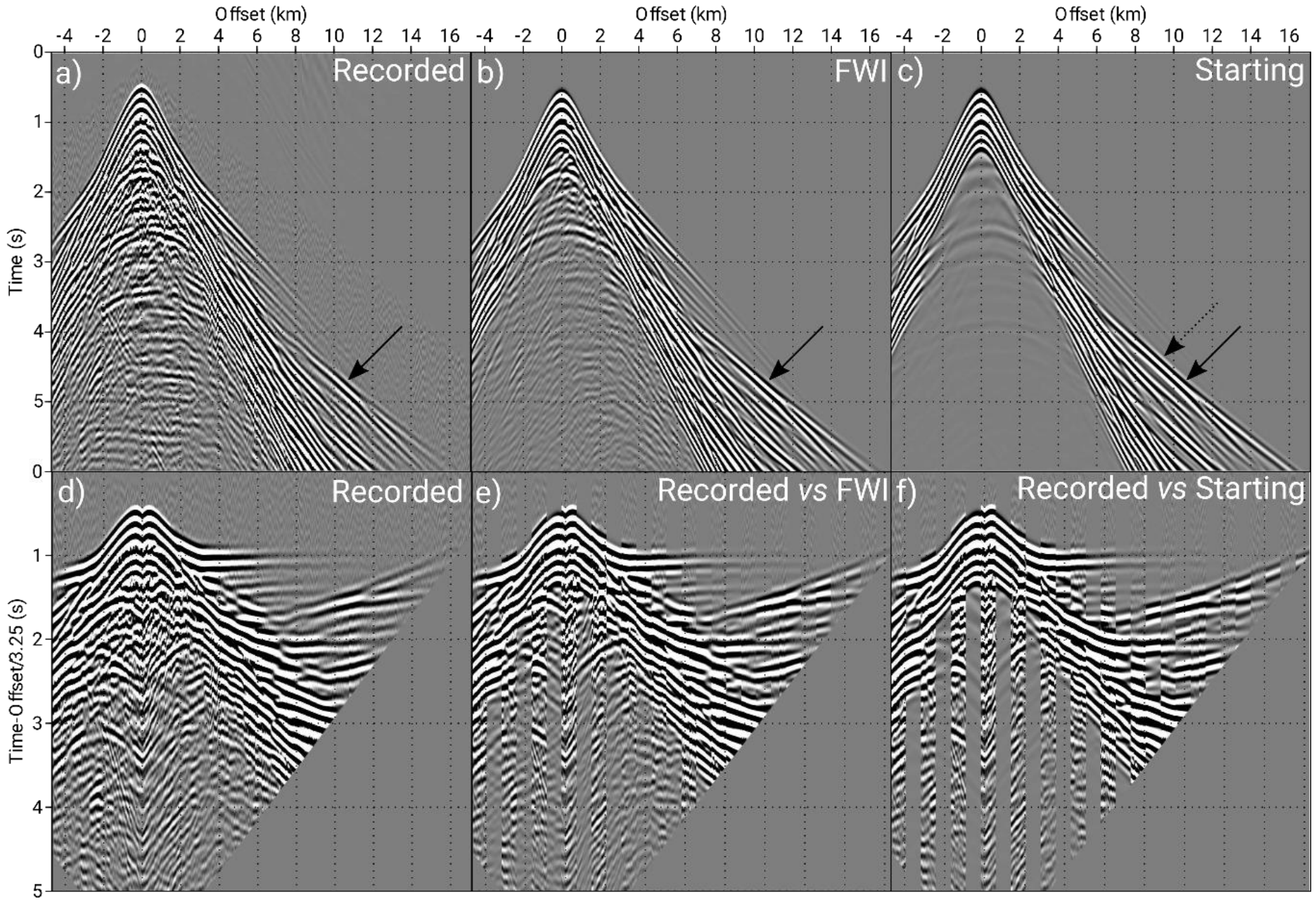

Figure 11: Time-domain recorded and simulated data (OBN 299). (a) Recorded data. (b-c) Synthetics computed in (b) FWI and (c) initial models. (d) Recorded data plotted with a reduction velocity (linear move out). (e) Interleaved recorded seismograms and computed seismograms in FWI model. (f) Interleaved recorded seismograms and computed seismograms in initial model. Solid black arrow points approximately critical distance of the recorded reflection from top of the reservoir (max amplitude). The dashed arrow approximately points the critical distance of the simulated reflection in the initial model from top of reservoir.

## Computational cost of frequency-domain FWI

We performed FWI on the Jean-Zay supercomputer of the national Institut du développement et des ressources en informatique scientifique (IDRIS). The peak power of the supercomputer is 28 PFlop/s. The Central Processing Unit (CPU) partition contains 1528 nodes with two processors Intel Cascade Lake 6248 (20 cores at 2,5 GHz), hence 40 cores by node. The shared memory per node is 192 GigaBytes. We allocate one MPI process per node and use 40 threads per MPI process to mitigate memory overheads during LU factorization. The size of the grids, the number of computer nodes used, the total memory and elapsed time required by the LU factorization and the elapsed time to compute 650 wavefields by forward/backward elimination are reviewed in Table **1** for several frequencies between 2.5 Hz and 13 Hz and for the two targets. In this table, classical full-rank factorization is used while the sparsity of the RHS is already exploited during the solution phase. From this table, one can see that the memory and the time costs of the LU factorization ($M_{LU}$ and $T_{LU}$, respectively) are manageable on large-scale supercomputers while the solution step is extremely fast ( $T_s$) thanks to the efficient block-processing of the sparse RHSs. We illustrate the impact of BLR and Mixed Precision (MP) in Figure **12**, which shows that the BLR approach implemented with mixed precision arithmetic (MP-BLR) further reduces quite significantly the memory (up to ≈25.6%) and time (up to ≈40%) of the LU factorization. Also, the highest cost reduction by MP-BLR is obtained for the largest computational domain in terms of covered area (target 1, 8.55 Hz), which is consistent with the admissibility condition used to compress the LU factors. This trend is indeed beneficial to tackle large-scale problems.

Table 1: f(Hz): Frequency; h(m) & ndof($10^6$): Grid interval and number of points in finite-difference grid; #*nodes*: Number of computer nodes; $M_{LU}(Mb)$: Total memory for LU factorization. $T_{LU}(s)$: Elapsed time for one LU factorization; $T_s(s)$: Elapsed time to compute 650 wavefields.

| f(Hz) | h(m) | ndof | #*nodes* | $M_{LU}$ | $T_{LU}$ | $T_s$ |
|---|---|---|---|---|---|---|
| Target 1 (23.5 $km$ × 30 $km$ × 8 $km$) | | | | | | |
| 2.5 | 150 | 2.62 | 16 | 72589.88 | 17.89 | 11.96 |
| 3.5 | 100 | 7.64 | 18 | 330207.53 | 87.03 | 25.40 |
| 5 | 75 | 16.51 | 34 | 923939.38 | 245.00 | 70.32 |
| 6.7 | 56 | 36.91 | 60 | 2825358.25 | 706.95 | 134.66 |
| 7.6 | 50 | 50.52 | 80 | 4348014.50 | 1072.19 | 200.46 |
| 8.5 | 45 | 68.04 | 120 | 6675021.00 | 1593.56 | 351.51 |

| Target 2 (23.5 *km* × 11 *km* × 8 *km*) | | | | | | |
|---|---|---|---|---|---|---|
| 8.5 | 45 | 26.15 | 60 | 1983768.50 | 477.33 | 81.99 |
| 10.1 | 37.5 | 43.36 | 80 | 3879103.50 | 1046.95 | 142.30 |
| 11.6 | 32.5 | 64.59 | 110 | 6860576.50 | 1952.99 | 249.05 |
| 13 | 30 | 73.89 | 128 | 8335750.50 | 2217.90 | 234.75 |

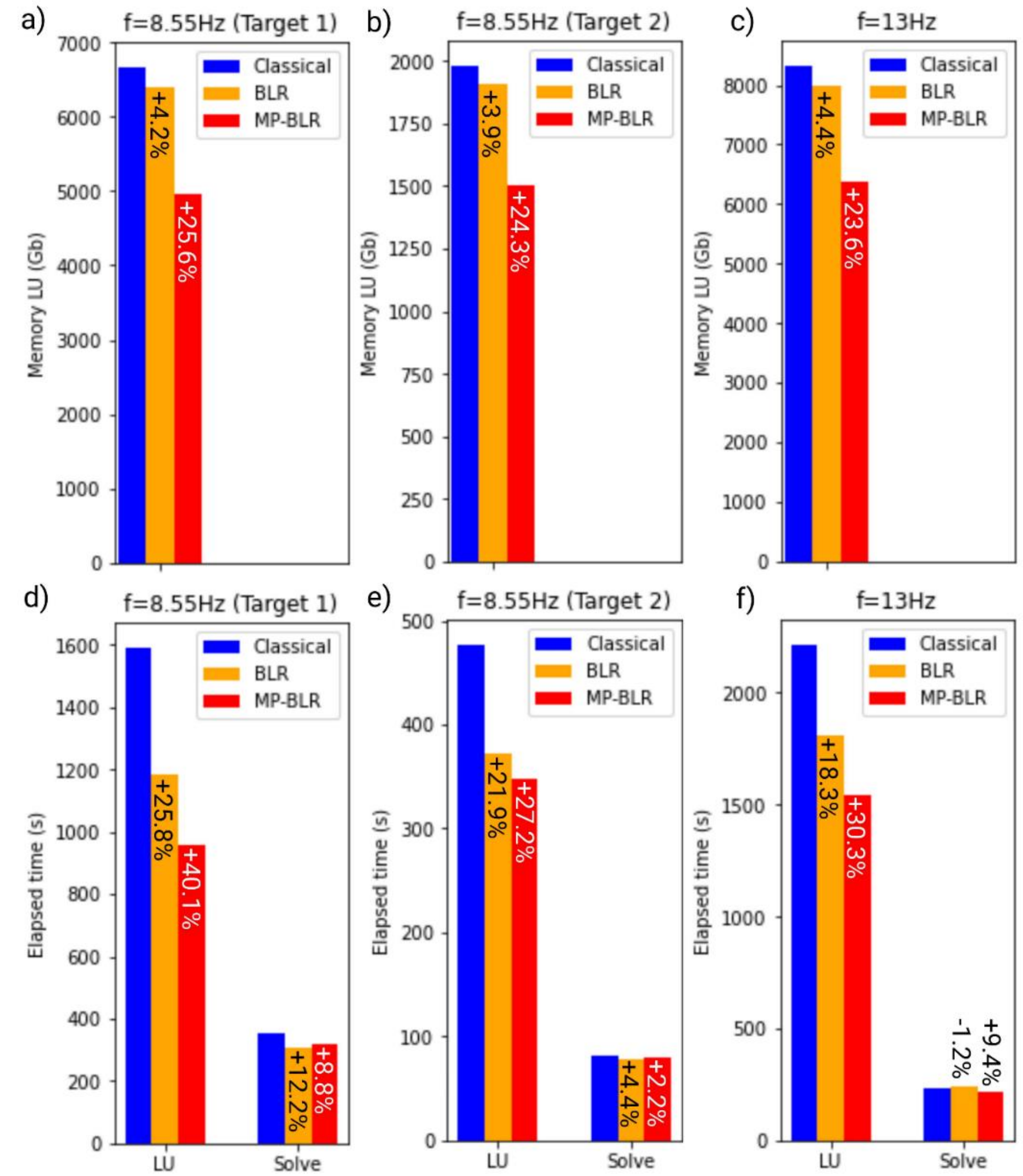

Figure 12: MUMPS statistics when performed with full-rank implementation, BLR and MP-BLR. (a-b) Memory cost of the LU factorization for frequencies (a) 8.55 Hz (Target 1); (b) 8.55 Hz (Target 2); (c) 13 Hz (Target 2). (d-f) Elapsed time for LU factorization and solution step for 650 RHSs. (d) 8.5 Hz frequency (Target 1). (e) 8.5 Hz frequency (Target 2). (f) 13 Hz frequency (Target 2). Note that the vertical scales are not the same in a,b,c and in d,e,f.

## CONCLUSION

We have shown the feasibility and the computational efficiency of 3D frequency-domain FWI based on direct solvers to tackle 3D OBN datasets at low to intermediate frequencies (< 15 Hz) with moderate computational resources (< 128 computer nodes). The strengths of the frequency-domain approach are two-fold: First, its versatility to design multiscale FWI with compact volume of data to rapidly build models with limited computational resources; Second, its ability to perform seismic modeling very efficiently for large number of sources. Moreover, attenuation can be efficiently implemented in FWI without computational overheads. We tackle FWI problems involving up to 74 million of parameters. However, we didn't reach the limits of our technology due to limited access to public computational resources. Therefore, it is likely that problems involving more than 100 million parameters can be tackled today with direct solvers that take advantage of low-rank compression and mixed precision arithmetic. Larger problems can be tackled with more scalable hybrid direct/iterative solvers based on domain decomposition preconditioner where a direct solver is used to solve local problems in each subdomain while the iterative solver is used to solve the global preconditioned system. However, the parallel processing of multiple sources will be less efficient with this approach. Aggressive frequency decimation, which is necessary for computational efficiency, may be detrimental to remove sparse acquisition footprint, which can be mitigated by sparsity-promoting regularization. We present results of mono-parameter FWI for $V_0$. Future works involve multi-parameter FWI for density and attenuation before considering visco-elastic FWI from the hydrophone and geophone components.

## ACKNOWLEDGMENTS

We thank Chevron and DUG for offering a copy of the pre-processed Gorgon OBN data and peripheral products to WIND. We are grateful to Chris Manuel (Chevron) for his assistance to download the Gorgon data. We are grateful to Gary Hampson (DUG) for his valuable comments that greatly help to improve the manuscript. This study was performed in the frame of the WIND project supported by AkerBP, ExxonMobil, Petrobras, Shell and Sinopec. The authors are grateful to the OPAL infrastructure from Observatoire de la Côte d'Azur for providing resources and support. This work was granted access to the HPC resources of GENCI under the allocation 0596.

## REFERENCES

Aghamiry, H., A. Gholami, L. Combe, and S. Operto, 2022, Accurate 3D frequency-domain seismic wave modeling with the wavelength-adaptive 27-point finite-difference stencil: a tool for full waveform inversion: Geophysics, **87**, 1–66.

Amestoy, P., O. Boiteau, A. Buttari, M. Gerest, F. J´ez´equel, J.-Y. L'Excellent, and T. Mary, 2022, Mixed-precision low-rank approximations and their application to block low-rank LU factorization: IMA Journal of Numerical Analysis, **https://doi.org/10.1093/imamum/drac037**, drac037.

Amestoy, P., R. Brossier, A. Buttari, J.-Y. L'Excellent, T. Mary, L. M´etivier, A. Miniussi, and S. Operto, 2016, Fast 3D frequency-domain FWI with a parallel Block Low-Rank multifrontal direct solver: application to OBC data from the North Sea: Geophysics, **81**, R363 – R383.

Amestoy, P., A. Buttari, J. Y. L'Excellent, and T. Mary, 2018a, On exploiting sparsity of multiple right-hand sides in sparse direct solvers: SIAM Journal on Scientific Computing, **41**, A269–A291.

———–, 2018b, On the complexity of the block low-rank multifrontal factorization: SIAM Journal on Scientific Computing, **49(4)**, A1710–A1740.

Bebendorf, M., 2004, Efficient inversion of the Galerkin matrix of general second-order elliptic operators with nonsmooth coefficients: Mathematics of computation, **251**, 1179–1199.

Brocher, T. M., 2005, Empirical relationships between elastic wavespeeds and density in the earth's crust: Bulletin of the Seismological Society of America, **95**, 2081–2092.

Duff, I. S., A. M. Erisman, and J. K. Reid, 2017, Direct methods for sparse matrices: Clarendon Press.

Gilbert, J. R., and J. W. Liu, 1993, Elimination structures for unsymmetric sparse LU factors: SIAM Journal on Matrix Analysis and Applications, **14**, 334–352.

Hicks, G. J., 2002, Arbitrary source and receiver positioning in finite-difference schemes using Kaiser windowed sinc functions: Geophysics, **67**, 156–166.

Métivier, L., and R. Brossier, 2016, The SEISCOPE optimization toolbox: A large-scale nonlinear optimization library based on reverse communication: Geophysics, **81**, F11– F25.

Métivier, L., R. Brossier, J. Virieux, and S. Operto, 2013, Full Waveform Inversion and the truncated Newton method: SIAM Journal On Scientific Computing, **35(2)**, B401–B437.

Operto, S., R. Brossier, L. Combe, L. M´etivier, A. Ribodetti, and J. Virieux, 2014, Computationally-efficient three-dimensional visco-acoustic finite-difference frequencydomain seismic modeling in vertical transversely isotropic media with sparse direct solver: Geophysics, **79(5)**, T257–T275.

Operto, S., and A. Miniussi, 2018, On the role of density and attenuation in 3D multiparameter visco-acoustic VTI frequency-domain FWI: an OBC case study from the North Sea: Geophysical Journal International, **213**, 2037–2059.

Operto, S., A. Miniussi, R. Brossier, L. Combe, L. Métivier, V. Monteiller, A. Ribodetti, and J. Virieux, 2015, Efficient 3-D frequency-domain mono-parameter full-waveform inversion of ocean-bottom cable data: application to Valhall in the visco-acoustic vertical transverse isotropic approximation: Geophysical Journal International, **202**, 1362–1391.

Plessix, R.-E., 2017, Some computational aspects of the time and frequency domain formulations of seismic waveform inversion, *in* Modern solvers for Helmholtz problems, Geosystems Mathematics: Springer, 159–187.

Saad, Y., 2003, Iterative Methods for Sparse Linear Systems: SIAM.

Sirgue, L., and R. G. Pratt, 2004, Efficient waveform inversion and imaging: a strategy for selecting temporal frequencies: Geophysics, **69**, 231–248.

Tournier, P.-H., P. Jolivet, V. Dolean, H. Aghamiry, S. Operto, and S. Riffo, 2022, Threedimensional finite-difference & finite-element frequency-domain wave simulation with multi-level optimized additive Schwarz domain-decomposition preconditioner: A tool for FWI of sparse node datasets: Geophysics, **87(5)**, T381–T402.

Virieux, J., and S. Operto, 2009, An overview of full waveform inversion in exploration geophysics: Geophysics, **74**, WCC1–WCC26.

Warner, M., A. Ratcliffe, T. Nangoo, J. Morgan, A. Umpleby, N. Shah, V. Vinje, I. Stekl, L. Guasch, C. Win, G. Conroy, and A. Bertrand, 2013, Anisotropic 3D full-waveform inversion: Geophysics, **78**, R59–R80.